\documentclass[draft]{article}
\usepackage[utf8]{inputenc}
\usepackage{amsmath,amssymb,amsfonts}
\usepackage{bm}
\usepackage{commath}

\usepackage{cleveref}

\usepackage{tikz,tikzscale}
\usepackage{graphicx}
\usepackage{pgfplots}
\usepgfplotslibrary{groupplots}
\pgfplotsset{compat=1.18}
\usepackage{pgfplotstable} 
\usetikzlibrary{backgrounds}
\usetikzlibrary{arrows,shapes}
\usetikzlibrary{mindmap}
\usetikzlibrary{chains}
\usetikzlibrary{calc}
\usetikzlibrary{matrix}
\usetikzlibrary{decorations.pathmorphing,patterns,shadows,decorations.markings}
\usepgfplotslibrary{fillbetween}
\usepgfplotslibrary{colormaps}
\makeatletter
\definecolor{gnuplot@orange}{RGB}{229,158,0}
\definecolor{gnuplot@purple}{RGB}{148,0,212}
\definecolor{gnuplot@lightblue}{RGB}{87,181,232}
\definecolor{gnuplot@green}{RGB}{0,158,75}
\definecolor{gnuplot@darkblue}{RGB}{0,115,179}
\definecolor{gnuplot@yellow}{RGB}{240,227,66}
\pgfplotscreateplotcyclelist{colorGPL}{%
gnuplot@darkblue,every mark/.append style={fill=gnuplot@darkblue!80!black},mark=o\\%
red!80!white,every mark/.append style={fill=red!80!black},mark=square\\%
gnuplot@green,every mark/.append style={fill=gnuplot@green!50!black},mark=otimes*\\%
gnuplot@orange,every mark/.append style={fill=gnuplot@orange!80!black,mark size=2.5pt},mark=x\\%
gnuplot@purple,mark=diamond\\%
black,densely dashed,every mark/.append style={solid,fill=gnuplot@darkblue!80!black},mark=square*\\%
gnuplot@lightblue,densely dashed,every mark/.append style={solid,fill=gnuplot@lightblue!30!black},mark=triangle*\\%
red!80!white,densely dashed,every mark/.append style={solid,fill=red!40!white},mark=oplus*\\%
}
\makeatother

\usepackage{Preamble/mypackages}
\usepackage{Preamble/mydefinitions}

\usepackage{ulem}

\title{Tensor-product vertex patch smoothers for biharmonic problems}
\author{Julius Witte\thanks{Interdisciplinary Center for Scientific Computing (IWR), Heidelberg University, Im Neuenheimer Feld 205, 69120 Heidelberg, Germany} \and Cu Cui\footnotemark[1]\hspace{4pt}\thanks{Supported by the China Scholarship Council (CSC) under grant no. 202106380059} \and Francesca Bonizzoni\thanks{MOX, Department of Mathematics, Politecnico di Milano, Piazza Leonardo da Vinci 32, 20133 Milano, Italy. FB has received support from the project PRIN2022, MUR, Italy, 2023--2025, P2022N5ZNP ``SIDDMs: shape-informed data-driven models for parametrized PDEs, with application to computational cardiology''.
FB is partially funded by  “INdAM - GNCS Project”, codice CUP E53C23001670001. FB is member of INdAM-GNCS} \and Guido Kanschat\footnotemark[1]\hspace{4pt}\thanks{Supported by the Deutsche Forschungsgemeinschaft (DFG) under Germany's Excellence Strategy EXC 2181/1 - 390900948 (the Heidelberg STRUCTURES Excellence Cluster)}}

\newcommand{\Ndim}{d}
\let\domain\Omega
\newcommand{\domainb}{{\partial\domain}}
\newcommand{\normal}{\bm{n}}
\newcommand{\tria}[1]{{\mathbb{T}_{#1}}}
\newcommand{\triah}{\tria{\ell}}
\newcommand{\dgfacesh}{\dgfaces{\ell}}
\newcommand{\dgfaceshi}{{\dgfaces{\ell}^{\circ}}}
\newcommand{\dgfaceshb}{\dgfaces{\ell}^{\partial}}

\newcommand{\dgspaceh}{\dgspace{\ell}} 

\DeclareMathOperator{\laplace}{\Delta}
\newcommand{\dbldot}{\mathbin{:}}
\newcommand{\dgmean}[1]{\left\{ #1 \right\}}
\newcommand{\forma}[2]{a(#1,#2)}
\newcommand{\formczeroip}[2]{a_{\ell}(#1,#2)} 
\newcommand{\seminormh}[1]{\envert{#1}_{\ell}}

\newcommand{\xx}{\bm{x}}

\newcommand{\dx}{\,\mathrm{d}x}
\newcommand{\dxx}{\,\mathrm{d}\xx}
\newcommand{\dsxx}{\,\mathrm{d}\sigma}

\pgfplotstableread{
fe	mv_2d	mv_3d	mv_2d_poisson	mv_3d_poisson	sm_2d	sm_2d_poisson	sm_3d	sm_3d_poisson
3	1.78E+09	1.03E+09	3.52E+09	3.82E+09	3.19E+08	6.05E+08	1.09E+08	3.75E+08
4	1.88E+09	1.17E+09	3.37E+09	3.09E+09	3.33E+08	5.83E+08	1.23E+08	3.02E+08
5	2.47E+09	1.15E+09	5.08E+09	2.32E+09	4.22E+08	7.89E+08	1.21E+08	2.34E+08
6	2.27E+09	9.89E+08	4.37E+09	2.05E+09	3.91E+08	6.89E+08	1.04E+08	2.13E+08
}\biharmop
\pgfplotstableread{
dof	fe	mv	smc1	gmresc1	sma1	gmresa1	smc2	gmresc2	sma2	gmresa2
1050625	2	2.25E+09	7.79E+08	3.78E+06	7.91E+08	3.90E+06	4.39E+08	5.57E+06	4.41E+08	5.70E+06
2362369	3	1.78E+09	6.72E+08	7.32E+06	6.76E+08	7.46E+06	4.05E+08	1.02E+07	4.08E+08	1.04E+07
4198401	4	1.88E+09	6.97E+08	1.04E+07	7.01E+08	1.03E+07	4.22E+08	1.50E+07	4.25E+08	1.48E+07
1640961	5	2.47E+09	8.17E+08	9.47E+06	8.29E+08	9.86E+06	4.68E+08	1.34E+07	4.74E+08	1.37E+07
2362369	6	2.27E+09	7.77E+08	9.83E+06	7.86E+08	1.01E+07	4.55E+08	1.22E+07	4.63E+08	1.26E+07
3214849	7	2.06E+09	7.25E+08	1.01E+07	7.30E+08	1.03E+07	4.31E+08	1.26E+07	4.34E+08	1.29E+07
}\avstwoD
\pgfplotstableread{
dof	fe	mv	sm	gmres	sm2	gmres2
1050625	2	2.25E+09	3.97E+08	8.47E+06	2.14E+08	5.92E+06
2362369	3	1.77E+09	3.19E+08	1.60E+07	1.75E+08	1.11E+07
4198401	4	1.88E+09	3.33E+08	1.77E+07	1.82E+08	1.22E+07
1640961	5	2.46E+09	4.22E+08	2.41E+07	2.29E+08	1.65E+07
2362369	6	2.27E+09	3.91E+08	2.29E+07	2.13E+08	1.64E+07
3214849	7	2.06E+09	3.56E+08	2.19E+07	1.95E+08	1.57E+07
}\mvstwoD
\pgfplotstableread{
dof	fe	mv	sma	gmresa	smc	gmresc	sma1	gmresa1
16974593	2	8.09E+08	2.33E+08	3.60E+06	2.31E+08	3.59E+06	3.61E+08	2.26E+06
57066625	3	1.03E+09	2.87E+08	5.43E+06	2.80E+08	5.31E+06	4.46E+08	4.17E+06
16974593	4	1.17E+09	3.06E+08	6.16E+06	3.02E+08	6.05E+06	4.84E+08	5.28E+06
33076161	5	1.15E+09	3.01E+08	5.51E+06	2.96E+08	5.46E+06	4.75E+08	5.05E+06
57066625	6	9.89E+08	2.62E+08	4.57E+06	2.60E+08	4.49E+06	4.12E+08	4.27E+06
90518849	7	8.24E+08	2.18E+08	4.00E+06	2.18E+08	4.00E+06	3.44E+08	3.25E+06
}\avsthreeD
\pgfplotstableread{
dof	fe	mv	sm	gmres	sm_e	gmres_e	sm2	gmres2
16974593	2	8.09E+08	8.70E+07	3.25E+06	8.06E+07	3.06E+06	4.58E+07	1.94E+06
57066625	3	1.03E+09	1.09E+08	6.44E+06	3.39E+07	2.07E+06	5.75E+07	4.89E+06
16974593	4	1.17E+09	1.23E+08	8.22E+06	1.00E+07	8.09E+05	6.46E+07	5.36E+06
33076161	5	1.15E+09	1.21E+08	7.94E+06	3.59E+06	3.72E+05	6.34E+07	6.42E+06
57066625	6	9.89E+08	1.04E+08	8.26E+06	1.68E+06	1.77E+05	5.48E+07	5.56E+06
90518849	7	8.24E+08	8.68E+07	6.99E+06	1.68E+06	1.77E+05	4.59E+07	4.72E+06
}\mvsthreeD
\pgfplotstableread{
dof2	perf2	perf2f	x2	dof3	perf3	perf3f	x3	dof4	perf4	perf4f	x4
35937	1.32E+06	1.56E+06	1.18E+00	117649	2.53E+06	3.17E+06	1.25E+00	274625	2.91E+06	3.73E+06	1.28E+00
274625	2.75E+06	4.14E+06	1.51E+00	912673	5.36E+06	8.12E+06	1.52E+00	2146689	6.47E+06	9.40E+06	1.45E+00
2146689	3.25E+06	5.34E+06	1.65E+00	7189057	5.87E+06	9.47E+06	1.61E+00	16974593	8.22E+06	1.29E+07	1.57E+00
16974593	3.28E+06	5.54E+06	1.69E+00	57066625	6.42E+06	1.06E+07	1.65E+00	135005697	8.21E+06	1.30E+07	1.59E+00
135005697	3.24E+06	5.51E+06	1.70E+00	454756609	6.43E+06	1.05E+07	1.63E+00	0	0	0.00E+00	0.00E+00
}\mixedprecision

\begin{document}
\maketitle

\begin{abstract}
We discuss vertex patch smoothers as overlapping domain decomposition methods for fourth order elliptic partial differential equations. We show that they are numerically very efficient and yield high convergence rates. Furthermore, we discuss low rank tensor approximations for their efficient implementation. Our experiments demonstrate that the inexact local solver yields a method which converges fast and uniformly with respect to mesh refinement and polynomial degree. The multiplicative smoother shows superior performance in terms of solution efficiency, requiring fewer iterations in both two- and three-dimensional cases. Additionally, the solver infrastructure supports a mixed-precision approach, executing the multigrid preconditioner in single precision while performing the outer iteration in double precision, thereby increasing throughput by up to 70\%.
\end{abstract}

\section{Introduction}

Recent years have seen a major adaptation of finite element software to modern hardware with its high levels of parallelism and an increasing gap between computing and memory speed. The key to this progress is high computational intensity. In this article, we address this issue by proposing a highly optimized multigrid smoother for the biharmonic problem and its implementation on a recent GPU accelerator.

The standard for efficient implementation of finite element methods has improved considerably in recent years~\cite{KronbichlerKormann12,KronbichlerKormann19,MelenkGerdesSchwab01} by the shift from sparse matrices to a computation of the action of the bilinear form to a vector on the fly. The amount of memory transfers was reduced considerably at the cost of additional computations on the data vector. Using sum factorization, the additional cost could be kept in check and the speed-up for applying the finite element operator to a vector was significant.
The development on the side of solvers was lagging behind for some time, as multigrid smoothers did either not provide the same computational intensity measured in floating point operations per byte loaded from main memory, or were very expensive.

In order to solve the biharmonic problem, we focus on an implementation of the $C^0$ interior penalty method (C0IP) of~\cite{BrennerSung05} albeit there is preceding work on multigrid for the biharmonic problem in~\cite{Zhang89,Zhao05}.
A multigrid method for C0IP with uniform convergence was developed in~\cite{BrennerZhao05}. In~\cite{BrennerWang05} two-level additive Schwarz preconditioners with a generic overlap at the scale of fine-level elements were analyzed. Similarly, in~\cite{FengKarakashian05}, theory on two-level non-overlapping additive and multiplicative Schwarz smoothers was provided, proving a bound for the condition number of $\bigo(H^3/h^3)$. The same bounds were demonstrated for isogeometric discretizations in~\cite{ChoPavarinoScacchi18}.

These works mainly focused on mathematical efficiency and did not address the question computational efficiency. To this end, we combine overlapping domain decomposition methods with multigrid. In order to have efficient local solvers, we choose vertex patches as subdomains, since on quadrilateral and hexahedral meshes they have tensor product structure. By~\cite{KanschatSharma14}, we can directly relate to our recent work on Stokes equations~\cite{cui2024stokes}. The key technique is the fast diagonalization method (FDM) from~\cite{LynchRiceThomas64}, which we apply in an approximate fashion as the local operators do not admit a rank two tensor representation.

Previous works~\cite{cui2024implementation,cui2024multilevel} have demonstrated the efficient implementation of vertex-patch smoothers on GPUs using the Poisson problem as a test case. By leveraging efficient shared memory utilization and optimized memory access patterns, these implementations achieved performance close to the shared memory roofline limits for most polynomial degrees. For instance, on an NVIDIA A100 GPU, the FP64 performance of vertex-patch smoothers and matrix-vector product kernels exceeded 3 TFLOPS, reaching 36\% of the theoretical peak performance, with similar efficiency observed for FP32. This approach was subsequently extended to solve the Stokes problem~\cite{cui2024stokes}.
Building on these advancements and taking advantage of the separable representation of the local solvers in the smoothers proposed in this work and its efficient evaluation with the fast diagonalization method, we follow the same implementation strategy to address biharmonic problems in this study.

In the following section, we introduce the model problem and its discretization. Section~\ref{sec:mg} describes the multigrid method and the additive and multiplicative vertex patch smoothers. In Section~\ref{sec:tensor_structure}, we discuss the tensor structure of the local problems which allows us to apply a fast approximate local solver. We present experimental results for various aspects of the method in Section~\ref{sec:experiments} before concluding.

\section{Model problem and discretization}
\label{sec:prelims}

In this article, we discuss a method for the model problem of the biharmonic equation
\begin{equation}
  \label{eq:model}
  \begin{aligned}
    - \laplace^2 u &= f &&\text{ in } \domain, \\
    u &= 0 &&\text{ on } \domainb, \\
    \nabla u \cdot \normal &= 0 &&\text{ on } \domainb, \\
  \end{aligned}
\end{equation}
where $\domain$ is a polygonal domain in $\mathbb{R}^\Ndim$ with
$\Ndim = 2, 3$.
The volumetric load $f$ is a given data. The boundary conditions are referred to as clamped.
The natural (weak) solution space for the biharmonic equation with clamped boundary conditions is $H^2_0(\domain)$.
The $C^0$-continuity and vanishing traces at the
physical boundary already implies that the tangential derivative vanishes as well at $\domainb$.
Then, the weak problem reads: find $u \in H^2_0(\domain)$ such
that
\begin{equation}
  \label{eq:weakform}
  \forma{u}{v} = F(v) \quad \forall v \in H^2_0(\domain),
\end{equation}
where
\begin{equation}
  \label{eq:bfandrhs}
  \forma{u}{v} = \int_\domain \nabla^2 u \dbldot \nabla^2 v \dxx,
  \qquad
  F(v)= \int_\domain f v \dxx.
\end{equation}
Here, the second order tensor $\nabla^2 v$ denotes the
Hessian matrix of a scalar field $v$ and the Hadamard product $A \dbldot B$ denotes
the full contraction of second order tensors $A$ and $B$.

Note that $a(u,u)$ is the square of the $L^2$-norm of the function $\nabla^2 u$. Applying Friedrichs' inequality twice, it thus introduces an inner product on $H^2_0(\Omega)$. Hence, the weak formulation~\eqref{eq:weakform} is well-posed.


The $C^0$ interior penalty method (C0IP) was introduced by Brenner and Sung~\cite{BrennerSung05, Brenner11} as a convenient discretization scheme for fourth order equations based on continuous finite element methods and a penalization of discontinuous first derivative by terms in the fashion of interior penalty methods~\cite{Arnold82,Nitsche71}.

We choose standard Lagrange finite elements on quadrilaterals and hexahedra with shape function spaces $\polyQ{\deg}$ consisting of (mapped) tensor product polynomials of degree $k$. Moreover, in order to prepare for multigrid methods, we discretize the domain $\domain$ by a nested sequence of meshes
\begin{gather*}
    \tria{0}\sqsubset\tria{1}
    \sqsubset\dots\sqsubset\tria{\ell}
    \sqsubset\dots\sqsubset\tria{L},
\end{gather*}
consisting of quadrilaterals or hexahedra of size $h_\ell$.
While our results are obtained on Cartesian meshes, the discretization only requires shape regular, locally uniform cells.
The relation ``$\sqsubset$'' indicates that all cells of a mesh $\tria{\ell}$ are obtained by refinement of the cells of $\tria{\ell-1}$. Refinement is obtained by dividing all edges of a cell in two and reconnecting in the natural way to $2^\Ndim$ new cells. On these meshes, we define finite element spaces $V_\ell$ by
\begin{equation}
  \label{eq:fespace}
  V_\ell = \left\{%
    v \in H^1_0(\domain) \;\middle|\;
    v_{|\cell} \in \polyQk
    \quad
    \forall \cell \in \tria{\ell}
  \right\}
  .
\end{equation}
From the definition it becomes apparent that $V_\ell$ is non-conforming, that is, $V_\ell$ is not a subspace of $H^2_0(\domain)$. Therefore, we follow~\cite{BrennerSung05} and introduce interface terms penalizing jumps in the derivative.

Let $\dgfaceshi$ and $\dgfaceshb$ denote the set of interior and boundary facets of $\tria{\ell}$, respectively, and let $\dgfacesh=\dgfaceshi\cup\dgfaceshb$. For an interior facet $e$, let $\cellp$ and $\cellm$ be the mesh cells such that $e = \cellp \cap \cellm$. On $e$, let $v^+$ and $v^-$ be the trace of a finite element function $v$ from these cells. Let $\partial_{n^+} v^+$ and $\partial_{n^-} v^-$ be the corresponding normal derivatives. Then, we define for $\xx\in e$ the jump of the normal derivative and the mean of the second derivative in normal direction as
\begin{gather}
\begin{split}
    \dgjump{\partial_n v}(\xx)
    &= \partial_{n^+} v^+(\xx) + \partial_{n^-} v^-(\xx),
    \\
    \dgmean{\partial_n^2 v}(\xx)
    &= \tfrac{1}{2} \left( 
    \partial^2_{n^+} v^+(\xx) + \partial^2_{n^-} v^-(\xx)
    \right).
\end{split}
\end{gather}
We emphasize that both definitions are independent of the choice of $\cellp$ and $\cellm$.
On a boundary facet, let
\begin{gather}
    \dgjump{\partial_n v}(\xx) = \partial_{n} v(\xx),
    \quad
    \dgmean{\partial_n^2 v}(\xx)
    = \partial^2_{n} v(\xx).
\end{gather}

By multiplying~\cref{eq:model} with a test function, integration by parts, and stabilizing by Nitsche's idea, the $C^0$ interior penalty ($C^0$-IP) formulation reads:
find $u_\ell \in \dgspaceh$ such that
\begin{equation}
  \label{eq:c0ip}
  \formczeroip{u_\ell}{v_\ell} = F(v_\ell) \quad \forall v_\ell \in V_\ell
\end{equation}
with the bilinear form
\begin{multline}
  \label{eq:bfc0ip}
    \formczeroip{u}{v} \defeq
     \sum_{K\in\tria{\ell}}\int_{K} \nabla^2 u \dbldot \nabla^2 v \dxx \\
     + \sum_{e\in\dgfacesh}\int_{e} \left(
      \frac{\ip}{h_\face} \, \dgjump{\partial_n u}  \dgjump{\partial_n v}
      - \dgmean{\partial_n^2 u}  \dgjump{\partial_n v}
      - \dgjump{\partial_n u}  \dgmean{\partial_n^2 v}
    \right) \dsxx
\end{multline}
and the linear form $F$ being defined in \eqref{eq:bfandrhs}.
 We refer to the cell integral of bilinear form~\cref{eq:bfc0ip} as
\emph{bulk} term and to the three face integrals, from left to right, as \emph{penalty}, \emph{consistency} and \emph{adjoint consistency} term. 

We assume that the penalty parameter $\ip$ only depends on the polynomial degree $k$ and $h_\face$ is the harmonic mean of the extent of the two mesh cells adjacent to $\face$ orthogonal to $\face$.

We recall the following results from~\cite{Brenner11,BrennerSung05}
\begin{proposition}[Well-posedness]
Assume that $\ip$ is chosen sufficiently large. Then, the bilinear form $\formczeroip{\cdot}{\cdot}$ is coercive with respect to the mesh-dependent norm
\begin{equation}
  \label{eq:hnorm}
  \normh{v}^2
  = \sum_{\cell \in \triah} \Hrseminorm{2}{v}{\cell}^2
  + \sum_{\face \in \dgfacesh} \frac{\ip}{h_\face} \Ltwonorm{\dgjump{\partial_n v}}{\face}^2,
\end{equation}
on $V_\ell$. Hence,~\eqref{eq:c0ip} has a unique solution.


Furthermore, if the weak solution $u$ to~\cref{eq:weakform} belongs to $\Hr{\deg+1}{\domain}$, then the interior penalty solution $u_\ell$
  in~\cref{eq:c0ip} satisfies the energy error estimate
  \begin{equation}
    \label{eq:energyerror}
    \seminormh{u - u_\ell}
    \leq C h_\ell^{\deg-1 } \Hrseminorm{k+1}{u}{\Omega}.
  \end{equation}
\end{proposition}

\section{Multigrid and vertex-patch smoothers}
\label{sec:mg}

The generic structure of a single $V$-Cycle multigrid step is summarized in Algorithm \ref{alg:gmg}. Note that pre- and post-smoothing can be one or more steps of the smoothers described below. The matrices $A_\ell$ are obtained from the bilinear form in equation~\eqref{eq:bfc0ip} by choosing a basis for the finite element spaces $V_\ell$ in the standard way. The choice of tensor product basis functions will be detailed in Section~\ref{sec:tensor_structure}.
\begin{algorithm}
\caption{V-Cycle on level $\ell$}
\label{alg:gmg}
\begin{algorithmic}[1]
  \Procedure{MG$_{\ell}$}{$x_{\ell},b_{\ell}$}
  \If{$\ell=0$}
  \State{$x_0 \gets A_0^{-1} b_0$} \Comment{coarse grid solver}
  \EndIf
\State $x_{\ell} \gets S_{\ell} (x_{\ell}, b_{\ell})$
\Comment{pre-smoothing}
\State $b_{\ell-1} \gets I^\downarrow_{\ell-1} \bigl(b_{\ell} - A_{\ell} x_{\ell}\bigr)$ \Comment{restriction}
\State $e_{\ell-1} \gets \text{MG}_{\ell-1}(0,b_{\ell-1})$ \Comment{recursion}
\State $x_{\ell} \gets x_{\ell} +  I^\uparrow_{\ell-1} e_{\ell-1}$ \Comment{prolongation}
\State $x_{\ell} \gets S_{\ell} (x_{\ell}, b_{\ell})$
\Comment{post-smoothing}
\State \Return{$x_{\ell}$}
\EndProcedure
\end{algorithmic}
\end{algorithm}
The prolongation $I^\uparrow_{\ell-1}\colon V_{\ell-1}\to V_\ell$ is the matrix of the standard embedding operator and the restriction $I^\downarrow_{\ell-1}$ is its transpose.

\begin{figure}[tp]
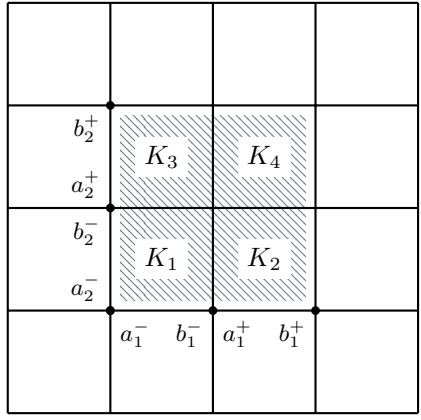
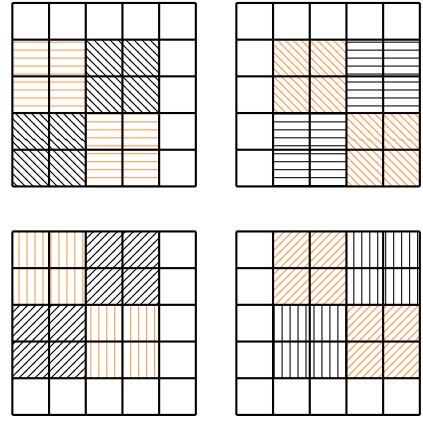

  \centering
  \begin{subfigure}{0.45\textwidth}
    \centering
    \includegraphics[width=\textwidth]{fig/left-right-cells-vp.tikz}
    \caption{Vertex patch $\vpatch=\{K_1, K_2, K_3, K_4\}$.}
  \end{subfigure}
  \hfill
  \begin{subfigure}{0.45\textwidth}
    \centering
    \includegraphics[width=.45\linewidth]{fig/coloring01-vp.tikz}
    \hspace{.05\textwidth}
    \includegraphics[width=.45\linewidth]{fig/coloring23-vp.tikz}
    \vspace{.02\textwidth}\\
    \includegraphics[width=.45\linewidth]{fig/coloring45-vp.tikz}
    \hspace{.05\textwidth}
    \includegraphics[width=.45\linewidth]{fig/coloring67-vp.tikz}
    \caption{Four nonoverlapping ``colors'' with additional red-black coloring for vertex patches.}
  \end{subfigure}
  \caption{A Cartesian vertex patch $\vpatch=\{K_1, K_2, K_3, K_4\}$ consisting of four cells in two dimensions (\emph{left}) and the coloring scheme (\emph{right}). The hatched and dotted areas on the right-hand
    shows a multiple red-black coloring satisfying $A_\ell$-orthogonality.}
  \label{fig:stencil-and-coloring-c0ip}
\end{figure}

The smoother is an overlapping domain decomposition method based on the multigrid method developed in~\cite{ArnoldFalkWinther97}. It is based on local solvers which solve the discretized differential equation on subdomains $\vpatch$ consisting of all cells attached to a vertex $\vertex$, so-called vertex patches, see also Figure~\ref{fig:stencil-and-coloring-c0ip}. By traversing all internal vertices of the mesh, the union of the patches covers the whole domain. This method has been applied successfully to Stokes and Brinkman equations~\cite{KanschatMao15,KanschatLazarovMao17} and it was shown that its application to the C0IP method is equivalent to the Stokes equations~\cite{KanschatSharma14}.

To be more precise, let $\vertex$ be an interior vertex of $\tria{\ell}$ and $\vpatch$ be its vertex patch, i.~e. the union of all cells adjacent to this vertex. Then, the associated subspace $V_{\vertex} \subset V_\ell$ consists of all finite element functions with support in $\vpatch$.
Here, we use the index $\ell$ for global vectors and omit it for local spaces and operators to keep the notation simple.
Let $R_{\vertex}$ be the localization operator, which selects out of a coefficient vector representing a finite element function in $V_\ell$ only those coefficients corresponding to functions in $V_{\vertex}$. Let furthermore $A_{\vertex}$ be the restriction of the matrix $A_\ell$ to $V_{\vertex}$. Later, we will also replace its inverse $A_{\vertex}^{-1}$ by an approximation. Then, we define the additive vertex patch smoother (AVS)
as
\begin{gather}
    S_{\ell;A} (x_{\ell}, b_{\ell})
    = x_{\ell} - \left[\sum_{\vertex} \omega R_{\vertex}^T A_{\vertex}^{-1} R_{\vertex}\right]
    \Bigl[b_\ell - A_\ell x_\ell \Bigr],
\end{gather}
where $\omega$ is a damping factor not exceeding the maximal number of patches overlapping in one cell, on uniform meshes 1/4 in two and 1/8 in three dimensions.
We will also provide examples, where the AVS is applied twice. Alternatively, we define the multiplicative vertex patch smoother (MVS). To this end, we define a numbering of the vertices $\vertex_1,\dots,\vertex_J$ and recursively
\begin{gather}
  \begin{split}
    x_\ell^{(0)} &= x_\ell,\\
    \nu=1,\dots,J\colon\quad x_{\ell}^{(\nu)}
    &= x_{\ell}^{(\nu-1)} - \omega R_{\vertex_\nu}^T A_{\vertex_\nu}^{-1} R_{\vertex_\nu}
    \bigl[b_\ell - A_\ell x_\ell^{(\nu-1)} \bigr],
    \\
    S_{\ell;M}(x_{\ell}, b_{\ell}) &= x_\ell^{(J)}. 
  \end{split}
\end{gather}
This algorithm, while not suffering from the necessity to choose a small damping factor, has two drawbacks. It is inherently sequential and there is a product with the matrix $A_\ell$ in every substep.
We fix the first issue by so-called coloring. This technique is inspired by the well-known red-black Gauss-Seidel method. First, we note that we obviously cannot write data of the same cell in a concurring way. Hence, in a first step, we split the set of all vertex patches into subsets of nonoverlapping patches called ``colors'', see the four meshes in Figure~\ref{fig:stencil-and-coloring-c0ip} on the right.
The computed update of each local vertex patch solver only affects the solution inside the patch. Hence, within one ``color'' these updates can be written at once and in parallel without fearing race conditions. But since the discontinuous Galerkin method reads from neighboring patches, we have to rely on that data to remain unchanged. This is achieved by the red-black coloring inside each nonoverlapping subset. Hence, on a uniform mesh, we obtain 8 colors in two and 16 colors in three dimensions.
Applying this scheme to the multiplicative algorithm with $n_c$ colors, and $\mathbb V_i$ being the set of all patches in color $i$, we obtain the following method:
\begin{gather}
  \begin{split}
    x_\ell^{(0)} &= x_\ell\\
    \nu=1,\dots,n_c\colon\quad x_{\ell}^{(\nu)}
    &= x_{\ell}^{(\nu-1)} - \left[\sum_{\vertex\in\mathbb V_\nu}\omega R_{\vertex_\nu}^T A_{\vertex_\nu}^{-1} R_{\vertex_\nu}\right]
    \bigl[b_\ell - A_\ell x_\ell^{(\nu-1)} \bigr],
    \\
    S_{\ell;M}(x_{\ell}, b_{\ell}) &= x_\ell^{(n_c)} 
  \end{split}
\end{gather}
Since the summations are nonoverlapping, we do not have the strict restriction on the damping factor $\omega$. Furthermore, we only have $n_c$ global operator applications instead of $J$. Finally, all the local solvers inside the sum can be applied in parallel. Hence, this is the implementation of the multiplicative smoother we use in the experiments.

\section{Tensor structure}
\label{sec:tensor_structure}

In this section we aim at highlighting the tensor product structure of the $C^0$-IP bilinear form~\eqref{eq:bfc0ip}, when restricted to a regular vertex patch. To this end, let us take rank-two trial and test functions $\varphi_\uindex{j}, \,\varphi_\uindex{i}$ of the form
\begin{equation}
  \label{eq:unitseparationofvariables}
  \begin{aligned}
  \varphi_\uindex{j}(x_1,x_2) &= \phi_{j_1}(x_1) \phi_{j_2}(x_2),\\    
  \varphi_\uindex{i}(x_1,x_2) &= \phi_{i_1}(x_1) \phi_{i_2}(x_2),\\
  \end{aligned}  
\end{equation}
where $\uindex{i}=(i_1,i_2)$ and $\uindex{j}=(j_1,j_2)$, and let us focus on the vertex patch $\vpatch\coloneqq\{K_1,\, K_2,\, K_3,\, K_4\}$ represented in Figure~\ref{fig:stencil-and-coloring-c0ip} (left).
Inserting $\varphi_\uindex{j},\,\varphi_\uindex{i}$ into~\eqref{eq:bfc0ip} restricted to $\mathcal{V}_K$, we get
\begin{align}
    \label{eq:bulk}
    &\sum_{K\in\vpatch} 
        \int_{K} \nabla^2 \varphi_\uindex{j} \dbldot \nabla^2 \varphi_\uindex{i} \dxx \,+\\
    \label{eq:boundary}
    & \sum_{K\in\vpatch}
    \sum_{e\in\partial K}\int_{e} \left(
      \frac{\ip}{h_\face} \, \dgjump{\partial_n \varphi_\uindex{j}}  \dgjump{\partial_n \varphi_\uindex{i}}
      - \dgmean{\partial_n^2 \varphi_\uindex{j}}  \dgjump{\partial_n \varphi_\uindex{i}}
      - \dgjump{\partial_n \varphi_\uindex{j}}  \dgmean{\partial_n^2 \varphi_\uindex{i}}
    \right) \dsxx.
\end{align}
Exploiting the tensor product form \eqref{eq:unitseparationofvariables} of the test and trial function, we can explicitly write the bulk term~\eqref{eq:bulk} as follows
\begin{align}
    \label{eq:bulk1}
  \eqref{eq:bulk}
  &=
  \left[\int_{a_1^-}^{b_1^-} \phi_{j_1}^{\prime\prime} \phi_{i_1}^{\prime\prime} \dx_1
  +\int_{a_1^+}^{b_1^+} \phi_{j_1}^{\prime\prime} \phi_{i_1}^{\prime\prime} \dx_1
  \right]
  \int_{a_2^-}^{b_2^+} \phi_{j_2} \phi_{i_2}\, \dx_2\\
  \label{eq:bulk2}
  &+2\int_{a_1^-}^{b_1^+} \phi_{j_1}^{\prime} \phi_{i_1}^{\prime} \dx_1
  \int_{a_2^-}^{b_2^+} \phi_{j_2}^{\prime} \phi_{i_2}^{\prime} \dx_2
  \\
  \label{eq:bulk3}
  &+\int_{a_1^-}^{b_1^+} \phi_{j_1} \phi_{i_1}\dx_1
  \left[
  \int_{a_2^-}^{b_2^-} \phi_{j_2}^{\prime\prime} \phi_{i_2}^{\prime\prime} \dx_2
  +\int_{a_2^+}^{b_2^+} \phi_{j_2}^{\prime\prime} \phi_{i_2}^{\prime\prime} \dx_2
  \right].
\end{align}
Analogously, we can reformulate the expression of the face integrals~\eqref{eq:boundary}. For instance, let $e$ be the vertical edge $e={b_1^-}\times [a_2^-,b_2^-]$, see also Figure~\ref{fig:stencil-and-coloring-c0ip} (left).
Then, the face integral 
\begin{gather}
\label{eq:e}
\int_{e} \left(
      \frac{\ip}{h_\face} \, \dgjump{\partial_n \varphi_\uindex{j}}
      \dgjump{\partial_n \varphi_\uindex{i}}
      - \dgmean{\partial_n^2 \varphi_\uindex{j}}
      \dgjump{\partial_n \varphi_\uindex{i}}
      - \dgjump{\partial_n \varphi_\uindex{j}}
      \dgmean{\partial_n^2 \varphi_\uindex{i}}
    \right) \dsxx
\end{gather}
becomes 
\begin{gather}
    \label{eq:ev}
    \left(\frac{\ip}{(|b_2^--a_2^-|)} \, \dgjump{\phi_{j_1}'}\dgjump{\phi_{i_1}'}
    -\dgmean{\phi_{j_1}''}\dgjump{\phi_{i_1}'}-\dgjump{\phi_{j_1}'}\dgmean{\phi_{i_1}''}
    \right) \int_{a_2^-}^{b_2^-} \phi_{j_2}\phi_{i_2}\dx_2. 
\end{gather}
For integrals over horizontal edges we have to swap the two dimensions. Namely, given $e= [a_1^-,b_1^-]\times b_2^-$, the face integral~\eqref{eq:e} becomes
\begin{gather}
    \label{eq:eh}
    \left(\frac{\ip}{(|b_1^--a_1^-|)} \, \dgjump{\phi_{j_2}'}\dgjump{\phi_{i_2}'}
    -\dgmean{\phi_{j_2}''}\dgjump{\phi_{i_2}'}-\dgjump{\phi_{j_2}'}\dgmean{\phi_{i_2}''}
    \right) \int_{a_1^-}^{b_1^-} \phi_{j_1}\phi_{i_1}\dx_1. 
\end{gather}
As a consequence, the  $C^0$-IP discretization matrix for a regular vertex patch has the following rank-3 tensor representation
\begin{equation}
  \label{eq:c0iptensorvp}
  B^\modeone \kp M^\modetwo
  +
  2 L^\modeone \kp L^\modetwo
  +
  M^\modeone \kp B^\modetwo,
\end{equation}
where $M^\moded$ and $L^\moded$ denote the mass and stiffness matrices along the direction $d$, namely, 
\begin{equation}
  \label{eq:matrix1d}
  M^\moded_{ij} = \int_{a_d^-}^{b_d^+} \phi_{j} \phi_{i}\, dx,
  \qquad
  L^\moded_{ij}
  =
  \int_{a_d^-}^{b_d^+} \phi_{j}^{\prime} \phi_{i}^{\prime}\, dx,
\end{equation}
and the matrices $B^\moded$ are computed from~\eqref{eq:bulk1}, \eqref{eq:bulk3}, \eqref{eq:ev} and~\eqref{eq:eh}.
In three dimensions, we obtain with the same arguments
\begin{multline}
      \label{eq:c0iptensorvp3D}
  B^\mode1 \kp M^\mode2 \kp M^\mode3
  + M^\mode1 \kp B^\mode2 \kp M^\mode3
  + M^\mode1 \kp M^\mode2 \kp B^\mode3\\
  + 2 L^\mode1 \kp L^\mode2 \kp M^\mode3
  + 2 M^\mode1 \kp L^\mode2 \kp L^\mode3
  + 2 L^\mode1 \kp M^\mode2 \kp L^\mode3.
\end{multline}

The above representation can also be utilized to perform matrix-vector multiplication. The matrix-free evaluation of the finite element operator $Ax$ using sum factorization has a computational complexity of $\mathcal{O}(d k^{d+1})$~\cite{KronbichlerKormann12,KronbichlerKormann19}. Since the finite element operator can be expressed as a sum of Kronecker products, as shown in~\eqref{eq:c0iptensorvp} and~\eqref{eq:c0iptensorvp3D}, this can be viewed as a specialized optimization. This approach involves fewer sum-factorization sweeps compared to standard numerical integration. 
It is important to note that these separable matrices are only applicable in the case of constant coefficients and axis-aligned meshes. The parallelization of this operator is implemented in a patch-wise manner~\cite{cui2024multilevel}, where partial cells and faces are processed within each patch. This approach not only avoids redundant integration computations but also improves data locality. Potential write conflicts can be avoided by so-called atomic operations. 

\subsection{Fast Diagonalization}
\label{sec:fastdiagonalization}

Here, we develop local solvers for the smoother which have the same complexity as the matrix-vector product above. To this end, we have to use the tensor structure again.
The fast diagonalization method was introduced in~\cite{LynchRiceThomas64} with the aim of efficiently solving second order separable partial differential equations in any dimension $d\geq 2$. For simplicity, let $d=2$. For any operator of the form
\begin{equation}
  \label{eq:laplacianfdm}
  A = A_1 \tp M_2 + M_1 \tp A_2,
\end{equation}
its inverse can be computed as
\begin{gather}
    \label{eq:inverse}
    A^{-1} = (Q_1\tp Q_2) \left(\Lambda_1\tp I_2 + I_1 \tp \Lambda_2\right)^{-1} (Q_1\tp Q_2)^\tpose,
\end{gather}
where the  diagonal matrices $\Lambda_i$ and the unitary matrices $Q_i$ solve the generalized eigenvalue problems
\begin{equation*}
Q_i^\tpose A_i Q_i = \Lambda_i M_i,\quad \text{for }i=1,\ldots,\Ndim.
\end{equation*}
Notice that in formula \eqref{eq:inverse} the inverse is taken over the diagonal matrices $\Lambda_i$ and the unitary matrices $Q_i$ can be applied separately in each direction. Hence, the application of $A^{-1}$ to a vector is indeed achieved with a complexity of $\mathcal O(dk^{d+1})$.

The operator \eqref{eq:c0iptensorvp} that represents the discretization of the biharmonic problem in $2\Ndim$ is not a rank-2 operator. Indeed, the elementary tensor $L^\modeone \kp L^\modetwo$ involving first order partial derivatives with respect to the first as well as second coordinate, prevents us from a separable tensor structure. A similar discussion holds in dimension $d=3$, where the terms of \eqref{eq:c0iptensorvp3D} that destroy the separable structure are $L^\mode1 \kp L^\mode2 \kp M^\mode3$, $M^\mode1 \kp L^\mode2 \kp L^\mode3$, $L^\mode1 \kp M^\mode2 \kp L^\mode3$. As a consequence, the exact local solvers are not amenable to fast diagonalization. 

On the other hand, if we omit the terms involving mixed first order partial derivatives arising from the bulk term, the local solvers on $\Omega_{\vertex_j}$ become inexact but with separable representation in two and three dimensions, respectively,
\begin{gather}
  \label{eq:localsolverbila}
  \begin{split}
    \tilde{A}_j &=B_j^\modeone \kp M_j^\modetwo+M_j^\modeone \kp B_j^\modetwo,
    \\
    \tilde{A}_j &=B_j^\modeone \kp M_j^\modetwo \kp M_j^\mode3 
    +M_j^\modeone \kp B_j^\modetwo \kp M_j^\mode3 
    + M_j^\modeone  \kp M_j^\modetwo \kp B_j^\mode3.  
  \end{split}
\end{gather}
Note that $B_j^\mode{i}$ and $M_j^\mode{i}$ are both symmetric as well as positive
semi-definite or positive definite, respectively, for $i=1,\dots,d$. 
Consequently, the generalized eigenvalue problems
required for fast diagonalization are well-defined. Besides, positive
definiteness of $\tilde{A}_j$ ensures that the local stability assumption is satisfied.

We close this section with a remark on non-Cartesian meshes. In this case, due to the mapping of vertex patches to a reference patch, even the Laplacian is not separable and the sums in equations~\eqref{eq:c0iptensorvp} and~\eqref{eq:c0iptensorvp3D} will not have the same low-rank structure anymore. Hence, the separable approximations in equation~\eqref{eq:localsolverbila} will become worse. This was documented for the Laplacian for instance in~\cite{WitteArndtKanschat20,Witte22}. A possible remedy is a change of the local basis functions as in~\cite{BrubeckFarrell22}. This discussion leads beyond the scope of this article and we defer it to a later publication.

\section{Computational experiments}
\label{sec:experiments}

We have implemented the numerical algorithms described above in the framework of the finite element library deal.II~\cite{dealII95}. The experiments use a system with an NVIDIA A100 GPU, hosted with two AMD EPYC 7282 16-core processors. The GPU has 80GB of high-speed HBM2e memory, providing 2TB/s peak memory bandwidth, and offers peak FP64 performance of 8.7 TFLOPS and FP32 performance of 17.4 TFLOPS on CUDA Cores and peak FP64 performance of 17.4 TFLOPS on Tensor Cores. 

Building on our previous work~\cite{cui2024implementation} on the Poisson problem, we will describe important practical details that were required to obtain performant implementations for these algorithms. 
To achieve high performance on GPUs, the most critical factor is the efficient utilization of the limited on-chip memory, commonly referred to as \emph{Shared Memory}, which offers lower latency and higher bandwidth. Following the approach from our previous work, we leverage shared memory to perform computations, thereby avoiding random high-latency global memory accesses.

The key to effectively utilizing shared memory lies in optimizing memory access patterns. By employing the conflict-free memory access pattern developed in our earlier work, optimal performance is achieved, where all memory load or store operations within a warp are serviced simultaneously without causing any \emph{bank conflicts}~\footnote{Shared memory is organized into multiple equally sized units called \emph{banks}, allowing simultaneous access by threads within a warp. If multiple threads within a warp access the same bank with different words, the accesses are serialized, resulting in a phenomenon known as \emph{bank conflict}.}. This optimal access pattern has been validated using the NVIDIA Nsight Compute profiler~\cite{cudaNcu}.
Finally, the entire solution process is performed exclusively on the GPU, eliminating data transfers between the CPU and the GPU, further enhancing performance and reducing latency.

When running a parallel for loop over all vertex patches of the additive smoother, the local operation on each patch can potentially lead to a race condition. This is because two overlapped patches will write or load the data from the same cell. Therefore, a similar colored algorithm used in the multiplicative algorithm can be applied to address this issue (referred to as \emph{colored AVS}). In this work, we also consider an alternative approach that leverages atomic operations to address this issue. On modern GPUs, atomic operations are executed with a single instruction using a ``fire-and-forget" semantic. This means that the instruction returns immediately, with conflict resolution managed by the cache system. Furthermore, modern GPUs provide support for fast atomic operations on double precision values. We denote this approach by \emph{atomic AVS}.

\begin{figure}[tbh]
  \centering
  \begin{tikzpicture}[scale=0.9]
    \begin{axis}[
      width=0.55\textwidth,
      xlabel={Polynomial degree $k$},
      ylabel={MDoF/s},
      tick label style={font=\footnotesize},
      label style={font=\footnotesize},
      legend to name=legendSM2d,
      legend columns=2,
      legend style={font=\footnotesize},
      title style={at={(1,0.952)},anchor=north east,draw=black,fill=white,font=\footnotesize},
      title={2D},
      xtick={2,3,4,5,6,7},
      ymin=0, 
      grid,
      cycle list name=colorGPL,
      mark size=1.8,
      semithick
      ]
      \addplot[draw=gnuplot@darkblue,mark=square] table[x={fe}, y expr=\thisrow{smc1}/1e6] {\avstwoD};
      \addlegendentry{colored AVS, one step};
      \addplot[draw=gnuplot@darkblue,densely dashed,mark=square,every mark/.append style={solid}] table[x={fe}, y expr=\thisrow{smc2}/1e6] {\avstwoD};
      \addlegendentry{colored AVS, two steps};
      \addplot[draw=red!80!white,mark=o,every mark/.append style={solid}] table[x={fe}, y expr=\thisrow{sma1}/1e6] {\avstwoD};
      \addlegendentry{atomic AVS, one step};
      \addplot[draw=red!80!white,densely dashed,mark=o,every mark/.append style={solid}] table[x={fe}, y expr=\thisrow{sma2}/1e6] {\avstwoD};
      \addlegendentry{atomic AVS, two steps};
      \addplot[draw=black,mark=otimes*,every mark/.append style={solid,fill=green!80!black}] table[x={fe}, y expr=\thisrow{sm}/1e6] {\mvstwoD};
      \addlegendentry{colored MVS, one step};
      \addplot[draw=black,densely dashed,mark=otimes*,every mark/.append style={solid,fill=green!80!black}] table[x={fe}, y expr=\thisrow{sm2}/1e6] {\mvstwoD};
      \addlegendentry{colored MVS, two step};
    \end{axis}
  \end{tikzpicture}
  \
  \begin{tikzpicture}[scale=0.9]
    \begin{axis}[
      width=0.55\textwidth,
      xlabel={Polynomial degree $k$},
      ylabel={MDoF/s},
      tick label style={font=\footnotesize},
      label style={font=\footnotesize},
      legend to name=legendSM3d,
      legend columns=2,
      legend style={font=\footnotesize},
      title style={at={(1,0.952)},anchor=north east,draw=black,fill=white,font=\footnotesize},
      title={3D},
      xtick={2,3,4,5,6,7},
      ymin=0, 
      grid,
      cycle list name=colorGPL,
      mark size=1.8,
      semithick
      ]
      \addplot[draw=red!80!white,mark=o,every mark/.append style={solid}] table[x={fe}, y expr=\thisrow{sma1}/1e6] {\avsthreeD};
      \addplot[draw=red!80!white,densely dashed,mark=o,every mark/.append style={solid}] table[x={fe}, y expr=\thisrow{sma}/1e6] {\avsthreeD};
      \addplot[draw=black,mark=otimes*,every mark/.append style={solid,fill=green!80!black}] table[x={fe}, y expr=\thisrow{sm}/1e6] {\mvsthreeD};
      \addplot[draw=black,densely dashed,mark=otimes*,every mark/.append style={solid,fill=green!80!black}] table[x={fe}, y expr=\thisrow{sm2}/1e6] {\mvsthreeD};
    \end{axis}
  \end{tikzpicture}
  \\
  \ref{legendSM2d}
  \caption{Throughput of the AVS and MVS smoothers with inexact local solver in two and three dimensions.}
\label{fig:smoother2d_3d}
\end{figure}
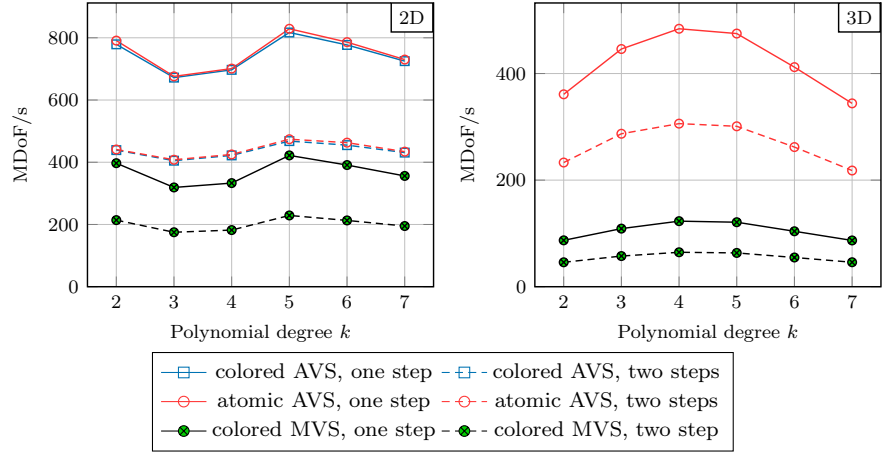

\Cref{fig:smoother2d_3d} compares the performance of different smoother implementations measured in MDoF/s. It can be observed that atomic AVS achieves performance comparable to colored AVS, reaching up to 800 MDoF/s for a single smoothing operation in two dimensions. In contrast, the more complex MVS is significantly affected by the additional residual updates within each color, with performance dropping notably in three dimensions to 100 MDoF/s, approximately one-third of that of AVS. Overall, for different polynomial degrees, the shared memory computation model employed ensures stable performance in both two- and three-dimensional cases.

\subsection{Numerical results}

We present results on Cartesian grids in two and three dimensions with the coarse grid being the decomposition of the square or cube $[0, 1]^d$ into $2^d$ cells, consequently, it consists of one vertex patch. We use the multigrid V-cycle as a preconditioner in the conjugate gradient solver (CG) and in the generalized minimal residual method (GMRES) for the additive and multiplicative versions, respectively. The stopping criterion for the Krylov subspace methods is a relative residual reduction of $10^{-8}$. On the coarse grid, we apply one smoothing step as the coarse grid solver.
The right-hand side $f$ is configured to yield an analytical solution of $u(x,y)=\sin(\pi x)\sin(\pi y)$ in two dimensions and $u(x,y,z)=\sin(\pi x)\sin(\pi y)\sin(\pi z)$ in three dimensions.
We measure the efficiency of the preconditioners by reducing the Euclidean norm of the residual $\|r_n\|$ after $n$ steps compared to the initial residual norm. Since the (integer) number of iteration steps depends strongly on the chosen stopping criterion,
we define the \emph{fractional iteration count} $\nu_{frac}$ by
\begin{equation*}
    \nu = -8 \log_{10} \bar{r}, \quad \bar{r} = \left(\frac{\|r_n\|}{\|r_0\|} \right)^{1/n}.
\end{equation*}

\begin{table}[tb]
  \caption{Fractional iterations $\fracsteps$ for additive 
    (AVS) and multiplicative (MVS) vertex patch smoothers in two dimensions, respectively with \emph{exact} local
    solvers. CG/GMRES solver to relative accuracy $10^{-8}$
    preconditioned by multigrid. Entries ``---'' not computed only levels $L$
    with $5 \times 10^4$ to $10^7$ DoFs.}
  \label{tab:biharm_cube_avs_and_mvs_exact}
  \centering
  \begin{tabular}{c ccccccc c}
    \toprule 
    Level $L$ & \multicolumn{6}{c}{Convergence steps $\fracsteps$} & Colors \\
    \cmidrule(lr){1-1} \cmidrule(lr){2-7} \cmidrule(lr){8-8}
    AVS (two steps) & $k=2$ & $k=3$ & $k=4$ & $k=5$ & $k=6$ & $k=7$ & \\
    \midrule 
       6  &   --- &   --- &  10.3 &  10.1 &  10.4 & 10.9 & 1 \\
       7  &  21.9 &  11.7 &   9.9 &   9.5 &  10.2 & 10.9 & 1 \\
       8  &  22.8 &  11.7 &   9.8 &   9.5 &   9.8 & 10.5 & 1 \\
       9  &  23.3 &  11.6 &   9.8 &   9.5 &   9.8 & 10.3 & 1 \\
      10  &  23.9 &  11.5 &   --- &   --- &   --- & --- & 1 \\
    \midrule
    MVS (one step) & $k=2$ & $k=3$ & $k=4$ & $k=5$ & $k=6$ & $k=7$ &  \\
    \midrule 
       6  &   --- &   --- &   2.9 &   2.6 &   2.5 & 2.3 & 8 \\
       7  &   8.8 &   4.4 &   2.9 &   2.6 &   2.5 & 2.3 & 8 \\
       8  &   8.9 &   4.4 &   2.9 &   2.5 &   2.4 & 1.9 & 8 \\
       9  &   9.2 &   4.4 &   2.9 &   2.5 &   2.4 & 1.9 & 8 \\
      10  &   9.3 &   4.4 &   --- &   --- &   --- & --- & 8 \\
    \bottomrule
  \end{tabular}
\end{table}

First, we show the numerical efficiency of the vertex patch smoother in terms of iteration steps without taking any tensor structure into account.
We perform two pre- and post-smoothing steps for the additive smoother with the natural relaxation parameter $ \omega = 1/4$ and one pre- and post-smoothing step for the multiplicative smoother with the relaxation parameter $\omega=1$. In~\Cref{tab:biharm_cube_avs_and_mvs_exact} iteration counts for both the
additive and multiplicative vertex patch smoother are compared. While both lead to uniform methods with respect to mesh level, the multiplicative method solves to high accuracy within 2--3 steps for polynomial degrees 4 and higher. Hence, they can be seen almost as direct solvers.


\begin{table}[tb]
  \caption{Fractional iterations $\fracsteps$ for additive 
    (AVS) and multiplicative (MVS) vertex patch smoothers in two dimensions, respectively with \emph{inexact} local
    solvers. CG/GMRES solver with relative accuracy $10^{-8}$
    preconditioned by multigrid. Entries ``---'' not computed only levels $L$
    with $5 \times 10^4$ to $10^7$ DoFs.}
  \label{tab:biharm_cube_avs_and_mvs_bila}
  \centering
  \begin{tabular}{c ccccccc c}
    \toprule 
    Level $L$ & \multicolumn{6}{c}{Convergence steps $\fracsteps$} & Colors \\
    \cmidrule(lr){1-1} \cmidrule(lr){2-7} \cmidrule(lr){8-8}
    AVS (two steps) & $k=2$ & $k=3$ & $k=4$ & $k=5$ & $k=6$ & $k=7$ & \\
    \midrule 
       6  &   --- &   --- &   9.2 &   9.9 &  10.3 & 10.8 & 1 \\
       7  &  19.2 &  10.4 &   9.0 &   9.3 &   9.8 & 10.3 & 1 \\
       8  &  19.6 &  10.3 &   9.0 &   9.2 &   9.5 & 9.9 & 1 \\
       9  &  19.9 &  10.3 &   9.0 &   9.1 &   9.4 & 9.8 & 1 \\
      10  &  20.1 &  10.2 &   --- &   --- &   --- & --- & 1 \\
    \midrule
    MVS (one step) & $k=2$ & $k=3$ & $k=4$ & $k=5$ & $k=6$ & $k=7$ & \\
    \midrule 
       6  &   --- &   --- &   4.2 &   4.5 &   5.2 & 5.8 &  8 \\
       7  &   9.2 &   4.8 &   4.2 &   4.5 &   5.2 & 5.8 &  8 \\
       8  &   9.4 &   4.8 &   4.2 &   4.5 &   5.1 & 5.7 &  8 \\
       9  &   9.5 &   4.8 &   4.2 &   4.4 &   4.9 & 5.7 &  8 \\
      10  &   9.5 &   4.8 &   --- &   --- &   --- & --- &  8 \\
    \bottomrule
  \end{tabular}
\end{table}

\begin{table}[!ht]
  \caption{Fractional iterations $\fracsteps$ for additive 
    (AVS) and multiplicative (MVS) vertex patch smoothers in three dimensions, respectively with \emph{inexact} local
    solvers. CG solver with relative accuracy $10^{-8}$
    preconditioned by multigrid. Entries ``---'' not computed only levels $L$
    with $4 \times 10^3$ to $10^7$ DoFs.}
  \label{tab:biharm_cube_avs_and_mvs_bila_3d}
  \centering
  \begin{tabular}{c ccccccc c}
    \toprule 
    Level $L$ & \multicolumn{6}{c}{Convergence steps $\fracsteps$} & Colors \\
    \cmidrule(lr){1-1} \cmidrule(lr){2-7}  \cmidrule(lr){8-8}
    AVS (one step) & $k=2$ & $k=3$ & $k=4$ & $k=5$ & $k=6$ & $k=7$ & \\
    \midrule 
       4  &   --- &   --- &   20.6 &  20.5  &   21.0 & 20.8 & 1  \\
       5  &   29.8 &  23.6 &   22.3 &   22.8 &  23.5 & 24.3 &  1  \\
       6  &  34.2 &  28.1 &   26.0 &   26.0 &   26.6 & 29.0 &  1  \\
       7  &  43.0 &  32.9 &   29.6 &   28.9 &   30.0 & 33.4 &  1  \\
       8  &  55.0 &  37.9 &   --- &   --- &   ---  & --- &  1 \\
    \midrule
    AVS (two steps) & $k=2$ & $k=3$ & $k=4$ & $k=5$ & $k=6$ & $k=7$ & \\
    \midrule 
       4  &   --- &   --- &   14.0 &   14.1 &   14.6 & 14.6 & 1  \\
       5  &   17.9 &  16.1 &   15.1 &   15.7 &  16.2 & 16.3 &  1  \\
       6  &  20.4 &  17.5 &   15.9 &   16.5 &   16.9 & 17.5 &  1  \\
       7  &  22.3 &  18.9 &   16.6 &   16.7 &   17.1 & 17.8 &  1  \\
       8  &  23.9 &  19.4 &   --- &   --- &   ---  & --- &  1 \\
    \midrule
    MVS (one step) & $k=2$ & $k=3$ & $k=4$ & $k=5$ & $k=6$ & $k=7$ & \\
    \midrule 
       4  &   --- &   --- &   5.4 &   5.9 &   5.9 & 6.8 & 16 \\
       5  &   9.1 &   4.9 &   4.6 &   4.9 &   4.7 & 4.9 & 16  \\
       6  &   9.3 &   5.3 &   4.3 &   4.3 &  4.0 & 4.4 & 16  \\
       7  &  9.5 &  5.2 &   4.1 &   4.2 &   3.9 & 4.0 & 16  \\
       8  &  9.5 &  5.2 &   --- &   --- &   ---  & --- & 16 \\
    \midrule
    MVS (two steps) & $k=2$ & $k=3$ & $k=4$ & $k=5$ & $k=6$ & $k=7$ & \\
    \midrule 
       4  &   --- &   --- &   3.2 &   3.6 &   3.7 & 4.0 & 16 \\
       5  &   7.2 &   4.0 &   3.3 &   3.0 &   2.9 & 2.9 & 16  \\
       6  &   8.0 &   3.9 &   3.3 &   2.9 &  2.7 & 2.7 & 16  \\
       7  &  8.4 &  3.8 &   3.3 &   2.9 &   2.6 & 2.7 & 16  \\
       8  &  8.6 &  3.8 &   --- &   --- &   ---  & --- & 16 \\
    \bottomrule
  \end{tabular}
\end{table}

Comparing the number of solver iterations in~\Cref{tab:biharm_cube_avs_and_mvs_bila} to those for exact local solvers
in~\Cref{tab:biharm_cube_avs_and_mvs_exact}, we even observe a slight advantage for the inexact method in case of additive vertex patches. For low polynomial degrees $\deg$ both multiplicative vertex patch methods compare at similar levels but with a growing gap in favor of exact local solvers if degrees increase: in particular, using inexact local solvers the number of iterations grows when increasing the polynomial degree while it decreases using exact local solvers. Nevertheless, the increase in iteration counts is moderate and absolute numbers are still low. 
The computational cost and memory intensity of these inexact but fast diagonalizable local solvers are significantly reduced, as reported in our previous work~\cite{cui2024implementation}. Therefore, Schwarz smoothers with inexact solvers based on fast diagonalization are expected to be superior in computational efficiency.

In the three-dimensional case, as shown in~\Cref{tab:biharm_cube_avs_and_mvs_bila_3d}, we compare fractional iterations for two vertex-patch smoothers (AVS and MVS) with either one or two smoothing steps. First, for AVS, introducing a second pre- and post-smoothing step significantly accelerates convergence, reducing the iteration count by approximately half. However, more iteration steps are required to achieve convergence compared to the two-dimensional case. This occurs because the local solver in three dimensions becomes more inexact due to neglecting additional terms, as illustrated in~\Cref{eq:localsolverbila}. 
For MVS, the method remains robust, yielding uniform convergence with respect to mesh refinement and polynomial degree. Furthermore, employing two smoothing steps reduces iteration counts even more effectively; specifically, for polynomial degrees $k > 4$, only three iterations are required to reach convergence.
For AVS, the relaxation parameter can be at most $1/8$ due to the overlap in the three-dimensional case. It must be reduced additionally because we use inexact local solvers. Our experimental results suggest that a relaxation parameter of $0.1$ provides good results.
In the case of MVS, we adopt a damping factor of $0.7$, selected based on empirical observations. Generally, MVS continues to outperform AVS, typically requiring approximately five iterations to converge, about one-quarter the iterations required by AVS.

\subsection{Performance}

\begin{figure}[tbh]
  \centering
  \begin{tikzpicture}[scale=0.9]
    \begin{axis}[
      width=0.55\textwidth,
      xlabel={Polynomial degree $k$},
      ylabel={MDoF/s},
      tick label style={font=\footnotesize},
      label style={font=\footnotesize},
      legend to name=legendSol2d,
      legend columns=2,
      legend style={font=\footnotesize},
      title style={at={(1,0.952)},anchor=north east,draw=black,fill=white,font=\footnotesize},
      title={2D},
      xtick={2,3,4,5,6,7},
      ymin=0, 
      grid,
      cycle list name=colorGPL,
      mark size=1.8,
      semithick
      ]
      \addplot[draw=gnuplot@darkblue,mark=square] table[x={fe}, y expr=\thisrow{gmresc1}/1e6] {\avstwoD};
      \addlegendentry{colored AVS, one step};
      \addplot[draw=gnuplot@darkblue,densely dashed,mark=square,every mark/.append style={solid}] table[x={fe}, y expr=\thisrow{gmresc2}/1e6] {\avstwoD};
      \addlegendentry{colored AVS, two steps};
      \addplot[draw=red!80!white,mark=o,every mark/.append style={solid}] table[x={fe}, y expr=\thisrow{gmresa1}/1e6] {\avstwoD};
      \addlegendentry{atomic AVS, one step};
      \addplot[draw=red!80!white,densely dashed,mark=o,every mark/.append style={solid}] table[x={fe}, y expr=\thisrow{gmresa2}/1e6] {\avstwoD};
      \addlegendentry{atomic AVS, two steps};
      \addplot[draw=black,mark=otimes*,every mark/.append style={solid,fill=green!80!black}] table[x={fe}, y expr=\thisrow{gmres}/1e6] {\mvstwoD};
      \addlegendentry{colored MVS, one step};
      \addplot[draw=black,densely dashed,mark=otimes*,every mark/.append style={solid,fill=green!80!black}] table[x={fe}, y expr=\thisrow{gmres2}/1e6] {\mvstwoD};
      \addlegendentry{colored MVS, two step};
    \end{axis}
  \end{tikzpicture}
  \
  \begin{tikzpicture}[scale=0.9]
    \begin{axis}[
      width=0.55\textwidth,
      xlabel={Polynomial degree $k$},
      ylabel={MDoF/s},
      tick label style={font=\footnotesize},
      label style={font=\footnotesize},
      legend to name=legendSol,
      legend columns=2,
      legend style={font=\footnotesize},
      title style={at={(1,0.952)},anchor=north east,draw=black,fill=white,font=\footnotesize},
      title={3D},
      xtick={2,3,4,5,6,7},
      ymin=0, 
      grid,
      cycle list name=colorGPL,
      mark size=1.8,
      semithick
      ]
      \addplot[draw=red!80!white,mark=o,every mark/.append style={solid}] table[x={fe}, y expr=\thisrow{gmresa1}/1e6] {\avsthreeD};
      \addplot[draw=red!80!white,densely dashed,mark=o,every mark/.append style={solid}] table[x={fe}, y expr=\thisrow{gmresa}/1e6] {\avsthreeD};
      \addplot[draw=black,mark=otimes*,every mark/.append style={solid,fill=green!80!black}] table[x={fe}, y expr=\thisrow{gmres}/1e6] {\mvsthreeD};
      \addplot[draw=black,densely dashed,mark=otimes*,every mark/.append style={solid,fill=green!80!black}] table[x={fe}, y expr=\thisrow{gmres2}/1e6] {\mvsthreeD};
    \end{axis}
  \end{tikzpicture}
  \\
  \ref{legendSol2d}
  \caption{Throughput measured as the number of DoF solved per second for the AVS and MVS with inexact local solver in two and three dimensions. CG/GMRES solver with relative accuracy $10^{-8}$ preconditioned by multigrid.}
\label{fig:solver2d_3d}
\end{figure}
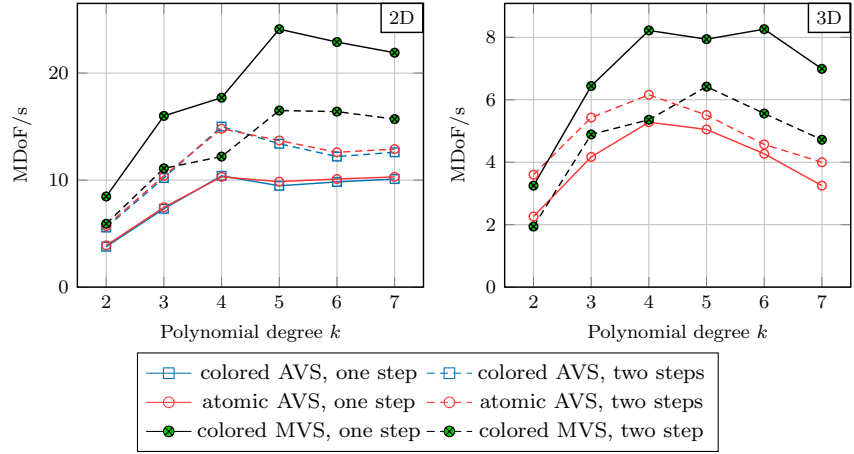

The solver throughput plotted in~\Cref{fig:solver2d_3d}, measured as the number of unknowns solved per second, shows that
the MVS is clearly superior to the AVS in two-dimensional case, especially for higher order elements. Since MVS requires fewer iterations, it demonstrates better performance in terms of solution efficiency, achieving more than 20 MDoF/s for polynomial degree $k>4$. For AVS, introducing a second pre- and post-smoothing step significantly enhances performance by reducing the required number of iterations. Similarly, AVS exhibits analogous behavior in three dimensions.

For MVS, superior and uniform performance continues to be observed in the three-dimensional case. However, although introducing a second smoothing step shows an advantage in terms of iteration counts (see~\Cref{tab:biharm_cube_avs_and_mvs_bila_3d}), the associated computational overhead, as illustrated in the right panel of~\Cref{fig:smoother2d_3d}, leads to a reduction in overall efficiency.

\begin{figure}[tbh]
  \centering
  \raisebox{-\height}{%
  \begin{tikzpicture}
    \begin{loglogaxis}[
      xlabel={Degrees of freedom},
      ylabel={MDoF/s},
      tick label style={font=\small},
      label style={font=\small},
      legend to name=legendMix,
      legend columns=1,
      legend style={font=\small},
      log basis y=2,
      ymin=1, ymax=20,
       ytick={1,2,4,8,16},
       yticklabels={1,2,4,8,16},
      grid,
      cycle list name=colorGPL,
      mark size=1.8,
      semithick
      ]
      \addplot[draw=gnuplot@darkblue,mark=square] table[x={dof2}, y expr=\thisrow{perf2f}/1e6] {\mixedprecision};
      \addlegendentry{$\mathbb{Q}_2$, mixed};
      \addplot[draw=gnuplot@darkblue,densely dashed,mark=square,every mark/.append style={solid}] table[x={dof2}, y expr=\thisrow{perf2}/1e6] {\mixedprecision};
      \addlegendentry{$\mathbb{Q}_2$, double};
      \addplot[draw=red!80!white,mark=o,every mark/.append style={solid}] table[x={dof3}, y expr=\thisrow{perf3f}/1e6] {\mixedprecision};
      \addlegendentry{$\mathbb{Q}_3$, mixed};
      \addplot[draw=red!80!white,densely dashed,mark=o,every mark/.append style={solid}] table[x={dof3}, y expr=\thisrow{perf3}/1e6] {\mixedprecision};
      \addlegendentry{$\mathbb{Q}_3$, double};
      \addplot[draw=black,mark=otimes*,every mark/.append style={solid,fill=green!80!black}] table[x={dof4}, y expr=\thisrow{perf4f}/1e6] {\mixedprecision};
      \addlegendentry{$\mathbb{Q}_4$, mixed};
      \addplot[draw=black,mark=otimes*,densely dashed,every mark/.append style={solid,fill=green!80!black}] table[x={dof4}, y expr=\thisrow{perf4}/1e6] {\mixedprecision};
      \addlegendentry{$\mathbb{Q}_4$, double};
    \end{loglogaxis}
  \end{tikzpicture}}  
  \raisebox{-\height}{\ref{legendMix}}
  \caption{Comparison of performance with double (solid) and mixed (dashed) precision for solving the biharmonic equation with GMRES preconditioned by multigrid with MVS with inexact local solver in three dimensions. Total solution time to relative accuracy $10^{-8}$.}
\label{fig:mixed}
\end{figure}
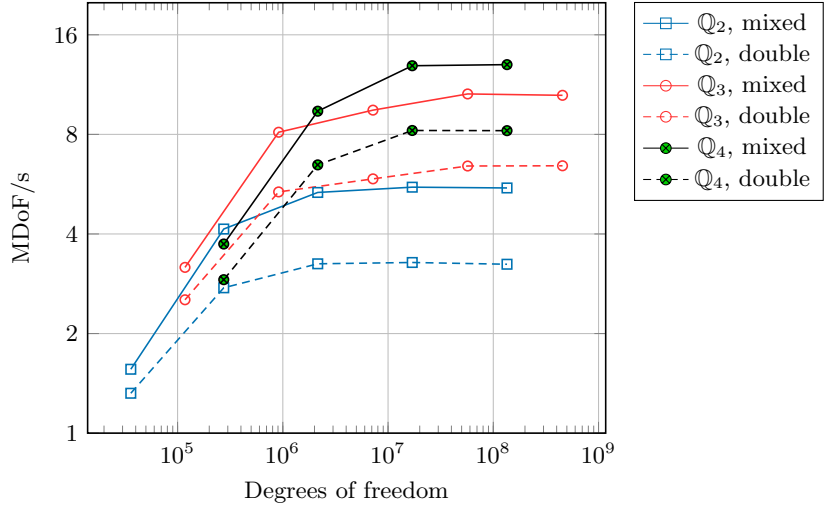

In iterative or direct solutions of linear systems, introducing lower precision is often employed to accelerate computations, benefiting from its high performance number. However, using single precision throughout the entire solution process is typically insufficient to achieve acceptable accuracy and a mixed-precision approach combining an outer double-precision CG iteration with a single-precision multigrid V-cycle as a preconditioner demonstrated a same level accuracy and iteration count of a fully double-precision reference~~\cite{goddeke2007performance}. Additionally, a recent work~\cite{CuiAcceleration} shown that using half-precision computations with NVIDIA Tensor Cores achieved over a 4x speedup.

In our implementation, the preconditioner, specifically the multigrid V-cycle, is executed entirely in single precision, with data converted at the entry point of the V-cycle. \Cref{fig:mixed} compares the performance of a GMRES iteration run in double precision, preconditioned by a multigrid V-cycle with MVS for pre- and post-smoothing in single precision, against performing all computations in the V-cycle using double precision.
For smaller problem sizes, both methods exhibit suboptimal performance as a result of insufficient parallelism to fully utilize the GPU. For larger problem sizes, we observed that the mixed-precision approach accelerated the solution of the linear system by 59\% to 70\%, capable of solving $8 \cdot 10^6$ unknowns per second in three dimensions.

\begin{figure}[tbh]
  \centering
  \begin{tikzpicture}[scale=0.9]
    \begin{axis}[
      width=0.55\textwidth,
      xlabel={Polynomial degree $k$},
      ylabel={MDoF/s},
      tick label style={font=\footnotesize},
      label style={font=\footnotesize},
      legend to name=legendPoisson,
      legend columns=2,
      legend style={font=\footnotesize},
      title style={at={(1,0.952)},anchor=north east,draw=black,fill=white,font=\footnotesize},
      title={$Ax$},
      ymin=0, 
      grid,
      cycle list name=colorGPL,
      mark size=1.8,
      semithick
      ]
      \addplot[draw=gnuplot@darkblue,mark=square] table[x={fe}, y expr=\thisrow{mv_3d}/1e6] {\biharmop};
      \addlegendentry{Biharmonic};
      \addplot[draw=black,mark=otimes*,every mark/.append style={solid,fill=green!80!black}] table[x={fe}, y expr=\thisrow{mv_3d_poisson}/1e6] {\biharmop};
      \addlegendentry{Poisson};
    \end{axis}
  \end{tikzpicture}
  \
  \begin{tikzpicture}[scale=0.9]
    \begin{axis}[
      width=0.55\textwidth,
      xlabel={Polynomial degree $k$},
      ylabel={MDoF/s},
      tick label style={font=\footnotesize},
      label style={font=\footnotesize},
      legend columns=2,
      legend style={font=\footnotesize},
      title style={at={(1,0.952)},anchor=north east,draw=black,fill=white,font=\footnotesize},
      title={MVS},
      ymin=0, 
      grid,
      cycle list name=colorGPL,
      mark size=1.8,
      semithick
      ]
      \addplot[draw=gnuplot@darkblue,mark=square] table[x={fe}, y expr=\thisrow{sm_3d}/1e6] {\biharmop};
      \addplot[draw=black,mark=otimes*,every mark/.append style={solid,fill=green!80!black}] table[x={fe}, y expr=\thisrow{sm_3d_poisson}/1e6] {\biharmop};
    \end{axis}
  \end{tikzpicture}
  \\
  \ref{legendPoisson}
  \caption{Comparison against Poisson problem for operator evaluation $Ax$ and one MVS step with inexact local solver in three dimensions with same level of refinement, i.e. on $64^3$ mesh.}
\label{fig:biharm_vs_poisson_3d}
\end{figure}
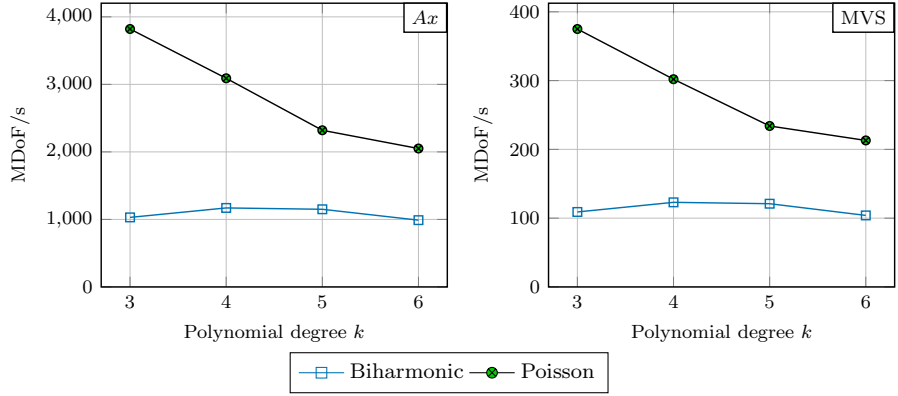

In our final evaluation of the efficiency of the GPU implementation presented in this work, we compare the performance of the finite element operator and the smoothing operator to that of the 3D Poisson problem, where the Poisson problem is discretized using the symmetric interior penalty method~\cite{cui2024multilevel}. We observed a 2–3x performance drop for the Biharmonic case to the Poisson problem. This is attributed to fact that the bi-Laplacian operator follows a rank-3 tensor representation~\eqref{eq:c0iptensorvp}, which involves additional terms, such as $L^\modeone \kp L^\modetwo$, compared to the Laplacian operator.
For the Poisson problem, the optimal performance is achieved at \(k=3\). This is because we use discontinuous Galerkin elements, which result in local operators and vector dimensions that are multiples of 8 (\(= 2k + 2\)), a configuration well-suited to GPU architectures. 

We believe that the implementation can be further optimized by better aligning with GPU hardware characteristics. For example, zero-padding could be applied during data reads and operator computations to enhance cache utilization and thread efficiency. However, it is important to note that this approach increases shared memory consumption, particularly in three-dimensional cases. Therefore, a careful balance must be struck to avoid excessive shared memory usage that could reduce occupancy and degrade operator performance.

\section{Conclusions}

We have considered the $C^0$-IP discretization of the biharmonic problem and we have presented exact vertex patch smoothers as well as an inexact version obtained by means of low rank tensor approximation. 
We have performed several numerical tests in two- and three-dimensional settings with the aim of showcasing the performances and demonstrating the efficiency of the proposed techniques. In particular, the approximation introduced by the low rank approximation still allows for a very fast converging multigrid methods with much less effort than the exact local solver. 
We observe that the multiplicative smoother demonstrates superior performance in both two- and three-dimensional cases. Meanwhile, for the additive smoother, introducing a second smoothing step significantly reduces the iteration count, thereby substantially enhancing the solver's performance.
Moreover, the computations have been executed in mixed precision: single precision for the multigrid preconditioner and double precision for the outer iteration of the multigrid preconditioner. The mixed-precision method solves the problem almost twice as fast.

\bibliographystyle{alpha}
\bibliography{references}
\end{document}